\documentclass[]{interact}
\usepackage{amssymb, amsthm, latexsym,  graphicx,makeidx, amsbsy,titlecaps, amsfonts,enumitem,mathtools}

\Addlcwords{of to with on a}

\usepackage{tikz}

\usepackage[caption=false]{subfig}% Support for small, `sub' figures and tables
\usepackage[natbibapa,nodoi]{apacite}
\usepackage[noabbrev]{cleveref}

\bibpunct[, ]{(}{)}{;}{a}{,}{,}% Citation support using natbib.sty
\renewcommand\bibfont{\fontsize{10}{12}\selectfont}% To set the list of references in 10 point font using natbib.sty

\newcommand{\sref}[1]{\subref{#1}}

\theoremstyle{plain}% Theorem-like structures provided by amsthm.sty
\newtheorem{theorem}{Theorem}[section]
\newtheorem{lemma}[theorem]{Lemma}

\theoremstyle{definition}

\theoremstyle{remark}

\newcommand{\N}{\ensuremath{\vec{\mathcal{N}}}}
\newcommand{\V}{\ensuremath{\Vert}}
\newcommand{\X}{\ensuremath{\mathbf{x}}}

\newcommand{\y}{\ensuremath{\mathbf{y}}}
\newcommand{\R}{\mathbb{R} }
\newcommand{\ddt}{ \frac{d}{dt}}

\newcommand{\abss}[1]{\V #1 \V}

\begin{document}

%\articletype{ARTICLE TEMPLATE}% Specify the article type or omit as appropriate

\title{\titlecap{An image segmentation model based on a variational formulation}}

\author{\name{Carlos M. Paniagua-Mej\'ia\textsuperscript{a}
%\thanks{} 
}\affil{\textsuperscript{a}Pontificia Universidad Cat\'olica Madre y Maestra, Dominican Republic; Email: cm.paniagua@ce.pucmm.edu.do} }
% \author{
% \name{Carlos M. Paniagua Mej\'ia\textsuperscript{a}
% \thanks{CONTACT A.~N. Author. Email: cm.paniagua@ce.pucmm.edu.do}}
% \affil{\textsuperscript{a}Pontificia Universidad Cat\'olica Madre y Maestra, Dominican Republic;Email: cm.paniagua@ce.pucmm.edu.do}

\maketitle

\begin{abstract}
Starting from a variational formulation, we present a model for image segmentation that employs both region statistics and edge information.  This combination allows for improved flexibility, making the proposed model suitable to process a wider class of images than purely region-based and edge-based models. We perform several simulations with real images that attest to the versatility of the model. We also show another set of experiments on images with certain pathologies that suggest opportunities for improvement.
\end{abstract}


\begin{amscode}53C44 35K55 49M25\end{amscode}

\section{Introduction}
In the last 30 years, beginning with the seminal works of  \cite{Kass88}, and  \cite{Shah89}, an ever-increasing number of variational and partial differential equation (PDE) based methods (often termed \textit{Deformable models} or \textit{Active Contours}) for image segmentation have been proposed. The basic idea is to overlay a contour over the given image and evolve it so that it stops at the boundaries of relevant objects present in the image. This is typically accomplished by minimizing some kind of \textit{energy} functional. Many other classical (non-PDE) approaches to the problem such as \textit{thresholding} and \textit{local filtering} (see e.g. \cite{Gonzalez04}) have also been introduced. However, the problem of segmenting an image is still open as there is no unifying theory within the field with methods cappable of attacking every image \cite[p.~367]{Marques11}. Our goal is to propose an image segmentation scheme that is able to segment two very important classes of images: piecewise constant images, and inhomogeneous images which appear quite frequently in the medical sciences \citep{Demirkaya09, Farag14}.

In what follows $\Omega $ represents a bounded and open subset of the plane $ \R^2 $. The grayscale image $ I \colon \Omega \rightarrow \R$ is realized as a bounded real valued function defined on $ \Omega $. In image processing tasks such as image segmentation ${\Omega}$ (along with its boundary) is typically a rectangle $[0,a]\times [0,b]$, and  $I$ is a discrete function taking values from 0 to 255.  The number $I(\X)=I(x,y)$  is called the \textit{graylevel} or intensity of $I$ at the  \textit{pixel} $\X = (x,y) \in \Omega$. By $C=C(q,t)\colon [0, 1]\times [0, \infty) \rightarrow \Omega$ we denote a smooth family of closed planar curves  where $q$ parameterizes the curve $C(\cdot,t)$ and $t$ parameterizes the family.

\section{Model Formulation}

    \subsection{Prior Models}\label{sec:preli}

\begin{figure}
\centering
\subfloat[Original image and initial curve.]{%
\resizebox*{5cm}{!}{\includegraphics{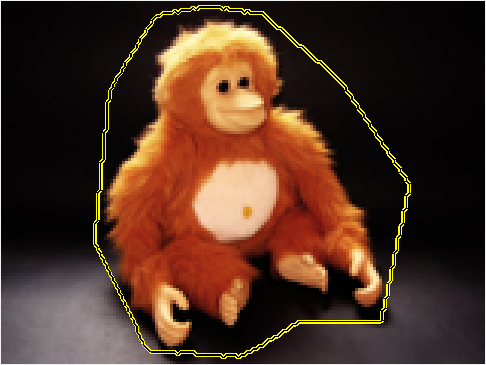}}}\hspace{5pt}
\subfloat[Final segmentation.]{%
\resizebox*{5cm}{!}{\includegraphics{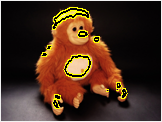}}}
\caption[Incorrect segmentation using Separation of Means.]{Incorrect segmentation using separation of means. The evolving curve is attracted to the brightest regions in the image.}
	\label{monkeymeansep}
\end{figure}

 To obtain a variational model for image segmentation a reasonable assumption about the composition of the image to be segmented is required.  In \cite{Yezzi} (see also \cite{Yezzi02}) it is assumed that $ I $ is \textit{bimodal}; that is, the image consists of two regions, foreground and background, with respective constant intensities $ \iota_1, \iota_2 $, with $ \iota_1\neq \iota_2.$
 %\footnote{We will relax this condition later.}
 Let $ \mu_1, \mu_2 $ be, respectively, the mean intensity of $ I $ inside and outside  $ C(\cdot,t) $. Then
 \begin{subequations}\label{eq:means1}
 	\begin{align}
 	\mu_1=\mu_1(\X,\Omega)&=\frac{1}{|\Omega|}\int_{\Omega} I(\X)~d\X,\label{means1}\\
 \mu_2=\mu_2(\X,\Omega)&=\frac{1}{|\Omega^c|}\int_{\Omega^c} I(\X)d\X, \label{means2}\end{align}
 \end{subequations} 
  \begin{equation}
 		|\Omega|=\int_\Omega d\X,
  \end{equation}
  where $ \Omega $ is the interior of $ C(\cdot,t) $, and $ \Omega^c $ its exterior. Let $ \omega $ represent the \textit{ideal} boundary separating the background and  foreground regions in $ I $
  (see \Cref{composition}).
  Clearly, as $ C $ approaches the unknown boundary $ \omega $ the absolute difference of the dynamic means $ \mu_1, \mu_2 $ increases and gets closer to the absolute difference of the region intensities $ \iota_1, \iota_2 $. \begin{figure}
\centering
	\begin{tikzpicture}
	\begin{scope}[blend group = soft light]
	\draw[blue!60!white,thick]   ( 90:1.2) circle (2);
	\fill[red!30!white]   ( 90:1.2) circle (2);
	\end{scope}
	\node at ( 90:2)    {Foreground};
	\node at ( 330:2)   {Background};
	\node[font=\Large] at (3,1.8) {$ \omega $};
	\draw[->, thick] (3,1.5) -- (2,1);
	 \draw [very thick, rounded corners, dashed] (-3,2) to[out=45,in=115] (-1,1) to[out=-180+115,in=10] (-3,-1) to [out=190,in=45+180](-3,2);
	 \node[font=\large] at (-4.2,1.5) {$C(\X,t)  $};
		\end{tikzpicture} \caption[Composition of a bimodal image.]{Composition of a bimodal image. Foreground (in color), background, and the boundary $ \omega $ (in blue) separating the two regions. A segmentation is obtained as $ C $  approaches $ \omega $.}
		\label{composition} \end{figure}
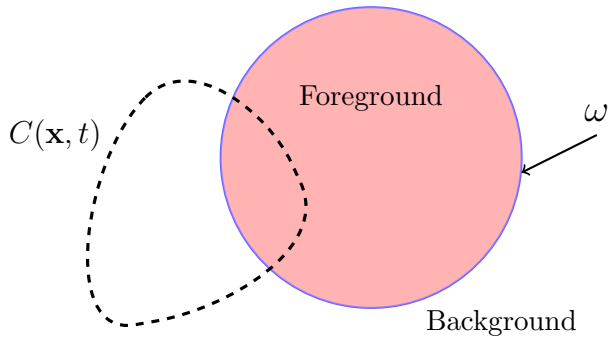 In symbols, \begin{equation}
	 |\mu_1 - \mu_2|  \rightarrow |\iota_1 - \iota_2| \text{ as } C \rightarrow \omega.
 \end{equation}Therefore, $ C $ gives a \textit{separation} of the foreground and background regions according to the `distance' separating their respective means. Based on this analysis, the following variational formulation is proposed (see \cite{Yezzi}): 
 
 \begin{equation}\label{mo:yezzi}
	\sup E= \sup ~(\mu_1-\mu_2)^2,
 \end{equation} where the $ \sup $ is taken over all admissible deformations of the curve $ C $ on $ W^{1,2}(\Omega)=\left\lbrace C\in L^2(\Omega) ;\nabla C \in L^2(\Omega)\right\rbrace $. We call this model the Mean Separation (MS) model as in \cite{Lankton08}. When the images are not bimodal, however, model \eqref{mo:yezzi} is not effective for segmentation of more general images (see \Cref{monkeymeansep}). 
 
    \subsection{Formulation}
  To extend the capabilities of model \eqref{mo:yezzi} to a larger class of images we  incorporate an \textit{edge function} $ g $   to augment its scope for boundary detection. We thus propose the following energy for object detection:
 \begin{equation}\label{edgeminsep1}
 	\hat{E}(\X,\Omega) = -\frac{1}{2}g(\X)(\mu_1-\mu_2)^2
 \end{equation} where $ \mu_1, \mu_2 $ are as in equations \eqref{means1} and \eqref{means2}. We expect the (possibly local) minima of \eqref{edgeminsep1} to attract evolving curves toward object boundaries. We say this model is a hybrid model because it uses two kinds of \textit{forces}--those exerted by the region features embedded in the \textit{separation} of the means $ \mu_1$ and  $\mu_2,$ and forces due to the edge function $ g $ that pins down the contours to salient discontinuities. To further strengthen \eqref{edgeminsep1} we add a regularization term (\cite{Tikhonov77})
 \begin{equation} \label{regula}
 	L=\oint_C ds
 \end{equation} 
 which is a boundary integral that exerts control over  $ C $ to have minimal arc length. We do so via a Lagrange multiplier $ \lambda>0 $ which controls the influence of the regularizing effects of \eqref{regula}. The resulting model is \begin{equation}\label{edgeminsep}
 		E(\X,\Omega) = -\frac{1}{2}~g(\X)~(\mu_1-\mu_2)^2+\lambda \oint_{\partial \Omega } ds.
 \end{equation}
 To obtain an evolution PDE we must compute the velocity field of \eqref{edgeminsep}. We do this next. 
 \subsubsection{Model Velocity Field}

  To implement our model for image segmentation we obtain the velocity field of \eqref{edgeminsep}. This can be done in either of two ways. One is the classical approach: convert region integrals to line integrals invoking solutions of Poisson's equation  with Dirichlet conditions with an application of Green's theorem. The other (more modern) alternative is to obtain the velocity field from the region functionals directly using shape derivation methods (see  \cite{Aubert03} and references therein). We proceed with the former approach and remark that the latter has interest in its own right. The region integrals in \eqref{edgeminsep} are of the form \begin{equation}\label{eq:regfunc}
 	F( \Omega) = \int_{\Omega}f(\X,\Omega)~ d\X
 \end{equation} where $ f $ is a scalar function and $ \Omega $ is an open, regular and bounded subset of $ \R^n $ with boundary $ \partial \Omega $. We prove the following result.

 \begin{lemma}\label{propo1}
 	Let $ \Omega $ be a bounded open set with regular boundary $ \partial \Omega $. 
 	Let $ f \colon \partial \Omega \rightarrow \R $ be continuous and u the unique solution of Poisson's equation:
	 \begin{equation}\label{eq:poisson}
	 \left\{
	 \begin{aligned}
	 -\Delta u &= f  \quad \text{in } \Omega,\\
	 u &= 0 \quad\text{on }  \partial \Omega.   
	 \end{aligned}\right.
	 \end{equation} The following identity holds: \begin{equation}\label{eq:converted}
	 	\int_{\Omega}f(\X,\Omega)~ d\X = \int_{\partial \Omega}\nabla u \cdot \N ~ds,
 	 \end{equation} where \N is the unit normal to $\partial \Omega.$
 \end{lemma}
 \noindent \textit{Proof.} If $ u $ is as in the hypothesis, then we have \begin{equation*}
	\int_{\Omega}f(\X,\Omega)~d\X  = -\int_{\Omega}\Delta u ~ d\X = \int_{\partial \Omega} \nabla u \cdot \N~ds.   
 \end{equation*} The first equality follows by \eqref{eq:poisson} inside $ \Omega $, and the second by Green's theorem. The proof is complete.
 \subsubsection*{Computation of the Gateaux derivative}
 With the region integrals in \eqref{edgeminsep} converted into line integrals, we next compute its Gateaux derivative, from which the velocity field will follow. We prove the more general result for functional \eqref{eq:regfunc}.
 \begin{theorem} \label{theorem1}
 	The Gateaux derivative of \eqref{eq:regfunc} is \begin{equation*}
 			F'( \Omega)  = -\int_{0}^{1} C_t \cdot f \N ~ds=-\int_{0}^{1} C_t(q) \cdot f(C(q))\N ~ds.
 	\end{equation*}
 \end{theorem} \noindent Note that $ F' $ does not depend on $ u $, but on the original integrand of functional \eqref{eq:regfunc}. 
 \noindent \\\textit{Proof}. By \Cref{propo1} we have
 \begin{equation}\label{fromproposition}
 		\int_{\Omega}f(\X,\Omega)~ d\X = \int_{\partial \Omega}\nabla u \cdot \N ~ds.
 \end{equation} Explicitly the gradient (column) vector\footnote{We use $\vec{V}^T=\langle a,b \rangle ^ T $ to denote the transpose of the row vector $ \vec{V} $. } field of $ u $ is \begin{equation}\label{ugradient}
 	\nabla u= \langle u_x,u_y \rangle^T,  
 \end{equation} and the inward unit normal vector $ \N $ to $ \partial\Omega $ is \begin{equation}\label{anormalvector}
 	\N=\frac{\langle -y_q,x_q \rangle^T }{\V C_q \V}.
 \end{equation} Substituting \eqref{ugradient} and \eqref{anormalvector} in \eqref{fromproposition} we get \begin{equation}
  \int_{\Omega}f(\X,\Omega)~ d\X = \int_{\partial \Omega}\nabla u \cdot \N ~dq  = \int_{0}^{1} (u_yx_q - u_x y_q) ~dq.
  \end{equation} Differentiating with respect to $ t $ gives \begin{equation}
  \begin{split}
  	F'( \Omega)  =\frac{d}{dt}	F( \Omega) &= \int_{0}^{1} \frac{d}{dt}(u_yx_q - u_x y_q) ~dq \\ &= \int_{0}^{1} (x_q\nabla u_y  \cdot C_t + u_y x_{qt} - ( y_q\nabla u_x \cdot C_t + u_x y_{qt})) ~dq \\
  	&= \int_{0}^{1} ((x_q\nabla u_y-y_q\nabla u_x)\cdot C_t +\langle u_y, - u_x\rangle\cdot C_{qt})~dq.
  \end{split}
  \end{equation} Integrating the second term of the last equation above by parts, one has \begin{equation}
  \begin{split}
  	F'( \Omega)  &=\int_{0}^{1} (x_q\nabla u_y-y_q\nabla u_x-\langle \nabla u_y \cdot C_q, -\nabla u_x \cdot C_q \rangle^T) \cdot C_t~dq \\
  	&=\int_{0}^{1}((u_{xx}+u_{yy})\langle -y_q, x_q \rangle^T)\cdot C_t ~dq\\
  	&=\int_{0}^{1} (\Delta u \N \V C_q \V)\cdot C_t~ dq.
  \end{split}
  \end{equation} The hypothesis of \Cref{propo1}, $ \Delta u = -f $, and the fact that $ ds=\V C_q \V dq $ verifies the claim. The proof is complete.

We now compute the Gateaux derivative of energy \eqref{edgeminsep}. Its first term is the region functional to which \Cref{theorem1} applies. The derivative of the second term is standard and given by
 \begin{equation}\label{eq:curveshort2}
 L'(t)=-\int_0^{L(t)} C_t \cdot \kappa \N ds,
 \end{equation} where $ \kappa $ is the mean curvature at each pixel \X ~of $ C $. Consider now the region dependent term of \eqref{edgeminsep}:
 \begin{equation}\label{eq:regfunc1}
 R(\X,\Omega)=-\frac{1}{2}~g(\X)~(\mu_1-\mu_2)^2,
 \end{equation} where the region functionals are embedded in the means $ \mu_1,\mu_2 $ of equation \eqref{eq:means1}. The Gateaux derivative of  \eqref{eq:regfunc1} is 
 \begin{equation}\label{derR}
  \begin{split}
  R'(\X,\Omega)&=-\frac{1}{2}~g(\X)~\frac{d}{dt}(\mu_1-\mu_2)^2\\
  &=g(\X)(\mu_2-\mu_1)\frac{d}{dt}(\mu_1-\mu_2)\\
  &=g(\X)(\mu_2-\mu_1)(\frac{d}{dt}\mu_1-\frac{d}{dt}\mu_2)
  \end{split} 
 \end{equation} where we have used the chain rule and the fact that $ g=g(I(\X)) $ does not change with time. The problem has been reduced to computing the derivatives of $ \mu_1 $ and $ \mu_2 $. We consider each individually beginning with $ \mu_1 $:
 \begin{equation}\label{derivationmu1}
 \begin{split}
  \ddt\mu_1 &=  \ddt\left(\frac{\int_{\Omega} Id\X}{\int_{\Omega} d\X}\right) \\
  &=\frac{\int_{\Omega} d\X ~\ddt \left( \int_{\Omega} Id\X \right) -\ddt\left( \int_{\Omega} d\X\right) \int_{\Omega} Id\X}{\left( \int_{\Omega} d\X\right)^2 }.
 \end{split}
 \end{equation} By \Cref{theorem1}, $\ddt \left( \int_{\Omega} Id\X \right) =-\int_0^1 C_t \cdot I\N ds $ and $ \ddt \left( \int_{\Omega} d\X \right)=-\int_0^1 C_t \cdot \N ds $. Also, to simplify notation, we write $ |\Omega|$ for $ \int_{\Omega} d\X $. Substituting these expressions in \eqref{derivationmu1} we get
 \begin{equation}\label{dermu1}
\begin{split}
 \ddt \mu_1 &=\frac{\int_0^1 C_t \cdot \N ds \int_{\Omega} I d\X-|\Omega|\int_0^1 C_t \cdot I\N ds}{|\Omega|^2}\\
&=\frac{\int_0^1 C_t \cdot \mu_1\N ds -\int_0^1 C_t \cdot I\N ds}{|\Omega|}\\
&=-\frac{\int_0^1 C_t \cdot (I-\mu_1)\N ds}{|\Omega|}\\
&=-\int_0^1 C_t \cdot \frac{I-\mu_1}{|\Omega|}\N ds
\end{split}
 \end{equation} For $ \mu_2=\int_{\Omega^c} Id\X~/\int_{\Omega^c}d\X $, we note that the inner normal vector to $ \Omega^c $ is $ -\N. $ Its Gateaux derivative is computed in similar fashion:
\begin{equation}\label{dermu2}
\ddt \mu_2 =\int_0^1 C_t \cdot \frac{I-\mu_2}{|\Omega^c|}\N ds.
\end{equation} Finally, using \eqref{dermu1} and \eqref{dermu2} in \eqref{derR} results (after some rearranging) in:
\begin{equation}\label{derRgood}
 R'(\X,\Omega)=-\int_0^1 C_t \cdot g(\X)(\mu_2-\mu_1)\left( \frac{I-\mu_1}{|\Omega|}+\frac{I-\mu_2}{|\Omega^c|} \right) \N ds,
\end{equation} and adding \eqref{eq:curveshort2} and \eqref{derRgood} we arrive at the Gateaux derivative for the proposed model:
\begin{equation}
	E'(\X,\Omega) =-\int_0^1 C_t \cdot \left[g(\X)(\mu_2-\mu_1)\left( \frac{I-\mu_1}{|\Omega|}+\frac{I-\mu_2}{|\Omega^c|} \right)+\lambda\kappa \right] \N ds.
\end{equation} 
% \subsubsection*{Construction of the velocity field for the solution of $ \min E $} 
From the Gateaux derivative we immediately obtain the associated Euler equation. Starting from an initial curve $C(\X,0)=C_0 $ the steepest descent method gives the following evolution equation for minima of \eqref{edgeminsep}:
\begin{equation}\label{edgeminsepevo}
C_t=g(\X)(\mu_2-\mu_1)\left( \frac{I-\mu_1}{|\Omega|}+\frac{I-\mu_2}{|\Omega^c|} \right)\N+\lambda\kappa  \N. 
\end{equation} We expect the steady state of \eqref{edgeminsepevo} to provide meaningful image partitions taking advantage of both edge and statistical information.

\section{\titlecap{Numerical implementation}}
We use the standard level set method of \cite{Osher88} for the numerical implementation of the model. We do not go into much detail for this derivation; we refer interested readers to \cite{Evans91} for a formal analysis and justification of the technique. Let $ C $ be represented by the zero level set of the embedding curve $ v \colon \R^2 \rightarrow \R $. Then
\begin{equation}\label{blah}
	\nabla v \cdot C_t + v_t=0.
\end{equation} Using \eqref{edgeminsepevo} in \eqref{blah} along with $ \N=-\nabla v/ \abss{\nabla v}  $ we obtain the evolution of $ C $ in terms of the level sets of $ v $:
\begin{equation*}
\begin{split}
v_t&=-	\nabla v \cdot C_t \\
&=	\nabla v \cdot \left( g(\X)(\mu_2-\mu_1)\left( \frac{I-\mu_1}{|\Omega|}+\frac{I-\mu_2}{|\Omega^c|} \right)+\lambda\kappa  \right) \frac{\nabla v}{\abss{\nabla v}}\\
&= \left( g(\X)(\mu_2-\mu_1)\left( \frac{I-\mu_1}{|\Omega|}+\frac{I-\mu_2}{|\Omega^c|} \right)+\lambda\kappa  \right)\abss{\nabla v}.
\end{split}
\end{equation*} In terms of $ v $ the curvature $\kappa$ of $C$ is given by $\kappa= \text{div}(\nabla v/ \abss{\nabla v}) $. Using this in the above equation we arrive at the level set evolution model for the solution of problem \eqref{edgeminsep}:

 \begin{equation}\label{theproblem2}
 \left\{
 \begin{aligned}
 &v_t=\left[ g(\X)(\mu_2-\mu_1)\left( \frac{I-\mu_1}{|\Omega|}+\frac{I-\mu_2}{|\Omega^c|} \right)+\lambda~ \text{div}\left(\frac{\nabla v}{\abss{\nabla v}} \right)   \right]\abss{\nabla v}\\
 &\frac{\partial v}{\partial \N} = 0 \text{ on } \partial \Omega \qquad \text{(Neumann boundary condition)}\\
 &v(\X,0)=v_0(\X) \qquad\text{(initial curve)} 
 \end{aligned}\right.
 \end{equation} \Cref{theproblem2} gives the evolution of all the levels sets of $ v $. However, we are only interested in the zero level set which represents $ C_t $. For our purposes it is sufficient to only consider pixels \X ~in a small neighborhood (strip or band) of the zero level set.  Aiming for computational efficiency, we implement localized techniques known as \textit{narrowband methods}. We refer the interested reader to references \cite{Sethian01,Sethian99,Adalsteinsson95} for details. For the discretization of \eqref{theproblem2} we employ forward differences for the time derivatives and central differences for the spacial derivatives at interior pixels and forward/backward differences for pixels at the left/right boundary accordingly. For the curvature $\kappa$ we have 
\begin{equation*}
\kappa= \text{div}\left( \frac{\nabla v}{\abss{\nabla v}} \right) =\frac{v_{xx} v_y^2-2v_{xy} v_x v_y + v_{yy} v_x^2}{(v_x^2+v_y^2)^{3/2}},
\end{equation*} so approximations to the second order partial derivatives $ v_{xx}, v_{yy}, $ and $ v_{xy} $ are needed for the approximation of $ \kappa $.  Let $ h $ represent the spatial step size on either direction. From Taylor's theorem we have the following expansions:\begin{equation}\label{expansion1}
v(x+h,y) = v(x,y)+hv_x(x,y)+\frac{h^2v_{xx}(x,y)}{2}+\frac{h^3v_{xxx}(x,y)}{6}+O(h^4)
\end{equation}
\begin{equation}\label{expansion2}
v(x-h,y) = v(x,y)-hv_x(x,y)+\frac{h^2v_{xx}(x,y)}{2}-\frac{h^3v_{xxx}(x,y)}{6}+O(h^4)
\end{equation} Adding \eqref{expansion1} and \eqref{expansion2} yields \begin{equation}
v(x+h,y)+v(x-h,y)=2v(x,y)+h^2 v_{xx}(x,y)+O(h^4)
\end{equation} and solving for $ v_{xx} $ gives the second order accurate formula \begin{equation}
v_{xx} (x,y)\approx\frac{v(x+h,y)-2v(x,y)+v(x-h,y)}{h^2}
\end{equation} which can be written using subindex notation as\begin{equation}
%\frac{v_{i+1,j}-2v_{ij}-v_{i,j-1}{h^2}
v_{xx}|_{ij}^n \approx \frac{v_{i+1,j}^n-2v_{ij}^n+v_{i-,j}^n}{h^2}
\end{equation} for every iteration $ n $ and interior pixel $ (ih,jh)=(i,j). $ Similarly for $ v_{yy} $ we get \begin{equation}
v_{yy}|_{ij}^n \approx \frac{v_{i,j+1}^n-2v_{ij}^n+v_{i,j-1}^n}{h^2}.
\end{equation} The second order mixed partial derivative $ v_{xy} $ requires more work. Refer to \Cref{fig:neighbors}. At the  four `corner' neighbors $ v(x+h,y+h), v(x-h,y-h), v(x-h,y+h),$ and $ v(x+h,y-h)$ we have the following expansions: 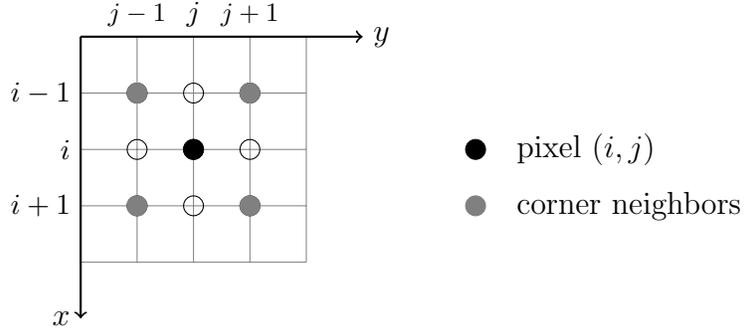
\begin{figure}
\centering
\begin{tikzpicture}[scale=.75]
\draw [help lines] (0,0) grid (4,-4);%grid
%axes
\draw[thick,->] (0,0)--(0,-5) node[anchor=east] {\large $x$};	%vertical axis
\draw[thick,->] (0,0)--(5,0) node[anchor=west] {\large $y$};	%horizontal axis
%interior pixels
\foreach \x in {1,2,3}
	\foreach \y in {-3,-2,-1}
		\draw  (\x,\y) circle (5pt);;
\filldraw (7,-2) circle (5pt) node[anchor=west] {\large \quad pixel $ (i,j) $};	
\filldraw[gray] (7,-3) circle (5pt); \node[right] at (7,-3) {\large \quad corner neighbors };	\filldraw[black] (2,-2) circle (5pt); 
\filldraw[gray] (1,-1) circle (5pt);
\filldraw[gray] (1,-3) circle (5pt);
\filldraw[gray] (3,-1) circle (5pt);
\filldraw[gray] (3,-3) circle (5pt);
%x labels
\node[anchor= east] at (0,-1) {$ i-1 $};
\node[anchor= east] at (0,-2) {$ i $};
\node[anchor= east] at (0,-3) {$ i+1 $};
%j labels
\node[anchor= south] at (1,0) {\small $ j-1 $};
\node[anchor= south] at (2,0) {$ j $};
\node[anchor= south] at (3,0) {\small $ j+1 $};
\end{tikzpicture} 
\caption[$ 3 \times 3 $ neighborhood of pixel $(i,j)  $.]{$ 3 \times 3 $ neighborhood of pixel $(i,j)  $. Corner neighbors are used for approximating mixed second order derivatives.}
\label{fig:neighbors}
\end{figure}

\begin{equation}\label{upright}
\begin{split}
-v(x+h,y+h) = &-(v( x,y)+h(v_x( x,y) +v_{y}(x,y)) \\
&+\frac{h^2 (v_{xx} ( x,y )+v_{yy}(x,y))}{2} +h^2 v_{xy}(x,y)) +\ldots
\end{split}
\end{equation}
\begin{equation}\label{downleft}
\begin{split}
-v(x-h,y-h) = &-(v( x,y)-h(v_x( x,y) +v_{y}(x,y)) \\
&+\frac{h^2( v_{xx} ( x,y )+ v_{yy}(x,y))}{2} +h^2 v_{xy}(x,y)) +\ldots
\end{split}
\end{equation}
\begin{equation}\label{downright}
\begin{split}
v(x+h,y-h) = \,&v(x,y)+h(v_x( x,y) -v_{y}(x,y)) \\
&+\frac{h^2 (v_{xx} ( x,y )+ v_{yy}(x,y))}{2} -h^2 v_{xy}(x,y) +\ldots
\end{split}
\end{equation}
\begin{equation}\label{upleft}
\begin{split}
v(x-h,y+h) = \,&v( x,y)-h(v_x( x,y) -v_{y}(x,y)) \\
&+\frac{h^2 (v_{xx} (x,y)+v_{yy}(x,y))}{2} -h^2 v_{xy}(x,y) +\ldots.
\end{split}
\end{equation} Adding equations \eqref{upright}--\eqref{upleft} and solving for $ v_{xy} $ yields the fourth order accurate formula
\begin{equation*}
v_{xy}=\frac{v(x+h,y+h)+v(x-h,y-h)-v(x+h,y-h)-v(x-h,y+h)}{4h^2}
\end{equation*} or more concisely using subindex notation
\begin{equation}
v_{xy}|_{ij}^n=\frac{v_{i+1,j+1}+v_{i-1,j-1}-v_{i+1,j-1}-v_{i-1,j+1}}{4h^2}.
\end{equation}

The edge function $g$ is chosen to be  \begin{equation} \label{edgestopping}
g = \frac{1}{1+\V \nabla G_{\sigma}\ast I \V},
 \end{equation}  
 where $G_\sigma \ast I$ is the convolution of the image $I$ with a Gaussian kernel  $G_\sigma =\sigma^{-1/2}e^{-|\X|^2/4\sigma}$ with standard deviation $\sigma.$
 \clearpage

 \section{Simulations}
  To assess the strengths and weaknesses of our variational model, which we call the \textit{Edge-Mean Separation} (EMS) model, in this section we demonstrate its performance on a number of images. In \Cref{simulations1} the examples illustrate how the EMS model is able to process images that the MS model is not able to segment correctly. Other technical aspects are tested in \Cref{simulations2}.

  \subsection{\titlecap{Application to wider class of images}}  \label{simulations1} 
\begin{figure}
        \centering
        \subfloat[]{%
        \resizebox*{.3\linewidth}{!}{\includegraphics{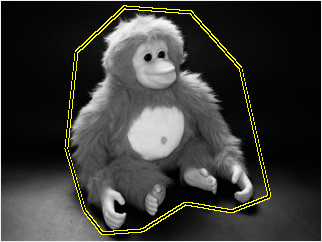}}}
        \subfloat[]{%
        \resizebox*{.3\linewidth}{!}{\includegraphics{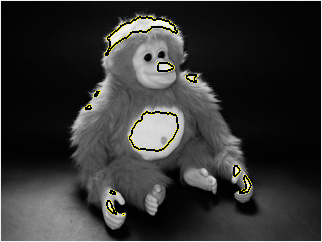}}}
        \subfloat[]{%
        \resizebox*{.3\linewidth}{!}{\includegraphics{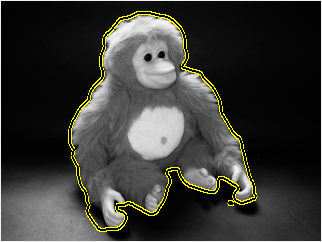}}}
        \caption[Successful segmentation of `monkey' image by our proposed model.]{Successful segmentation using our proposed model. (a) Image and initial curve; (b) result with the MS model; (c) result with the proposed EMS model.}
        	\label{monkeygoodseg}
        \end{figure}

Unless stated otherwise, we choose the model parameters $\lambda $ and  $\sigma$ to 1. We first revisit the `monkey' image for our first example. We noted (\Cref{sec:preli}) that the MS Model cannot delineate the boundaries of the relevant object for this image due to the presence of nonuniform lighting and bright regions within the object of interest. In \Cref{monkeygoodseg}, using the same initial curve as for the MS model, we obtain a meaningful segmentation of the monkey.

\begin{figure}
        \centering
        \subfloat[]{%
        \resizebox*{.23\linewidth}{!}{	\includegraphics[width=\textwidth]{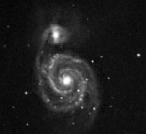}}}
        \subfloat[]{%
        \resizebox*{.23\linewidth}{!}{	\includegraphics[width=\textwidth]{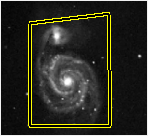}}}
        \subfloat[]{%
        \resizebox*{.23\linewidth}{!}{\includegraphics[width=\textwidth]{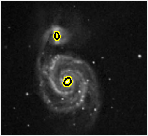}}}
        \subfloat[]{%
        \resizebox*{.23\linewidth}{!}{\includegraphics[width=\textwidth]{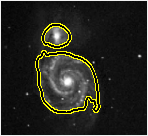}}}
        	\caption[A successful segmentation of `galaxy' image by our proposed model.]{A successful segmentation of `galaxy' image by our proposed model. The ``congnitive" boundary of the galaxy is suggested by our model. (a) Orignal image; (b) initial curve; (c) Result with MS model; (d) result with proposed EMS model.}
         	\label{im:galaxy}
        \end{figure}

 For our next experiment we consider a galaxy image. \Cref{im:galaxy} shows an example of a \textit{spiral} galaxy  with two spiral arms \citep{DeVaucouleurs59}. Our model is able to detect the `cognitive' boundary of this galaxy. Note further the presence of a large star at the end of the spiral arm extending toward the top of the image. Our model suggests the presence of a potential object of interest there along with nearby cosmic dust, gas, and other smaller stars--perhaps a second\textit{bulge} \citep{DeVaucouleurs59}.
 
\begin{figure}
        \centering
        \subfloat[]{%
        \resizebox*{.23\linewidth}{!}{	\includegraphics[width=\textwidth]{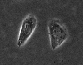}}}
        \subfloat[]{%
        \resizebox*{.23\linewidth}{!}{	\includegraphics[width=\textwidth]{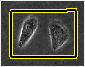}}}
        \subfloat[]{%
        \resizebox*{.23\linewidth}{!}{\includegraphics[width=\textwidth]{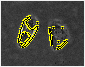}}}
        \subfloat[]{%
        \resizebox*{.23\linewidth}{!}{\includegraphics[width=\textwidth]{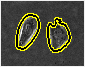}}}
        		\caption[Segmentation results for `two cells' image.]{(a) Original image and (b) initial curve; (c) results with the MS model and (d) the proposed EMS model. The two cells are successfully detected by our model.}
         \label{im:twocells}
        \end{figure}

  Our third and final example shows a medical image involving two cells, \Cref{im:twocells}. The advantages of using a level set formulation are apparent: starting from a rectangular front, a change in topology occurs and the front adjusts to the shape suggested by image's objects of interest. As in the previous examples our proposed model provides the more meaningful segmentation.
  
  \subsection{\titlecap{Further experiments}} \label{simulations2} In this section we show the results of several experiments conducted to assess other desired properties in segmentation techniques. Particularly, we investigate sensitivity to initialization, performance on \textit{multimodal} images, and sensitivity to noise.
  
  \subsubsection{Sensitivity to Initialization}

\begin{figure}
        \centering
        %a
        \subfloat[]{\label{correct1}
        \resizebox*{.45\linewidth}{!}{
            \includegraphics[width=.5\textwidth]{twocells_ini}%
            \includegraphics[width=.5\textwidth]{twocells_good}}}
        \subfloat[]{%b
        \resizebox*{.45\linewidth}{!}{	\includegraphics[width=.5\textwidth]{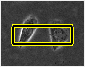}%    
         		\includegraphics[width=.5\textwidth]{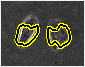}}}\\
        %c 		
        \subfloat[]{
        \resizebox*{.45\linewidth}{!}{	\includegraphics[width=.5\textwidth]{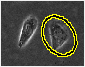}%    
         		\includegraphics[width=.5\textwidth]{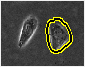}}}
         		\subfloat[]{
        \resizebox*{.45\linewidth}{!}{	\includegraphics[width=.5\textwidth]{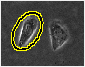}%    
         		\includegraphics[width=.5\textwidth]{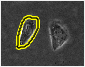}}}\\
         		\subfloat[]{\label{correct2}
        \resizebox*{.45\linewidth}{!}{	\includegraphics[width=.5\textwidth]{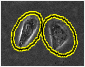}%    
         		\includegraphics[width=.5\textwidth]{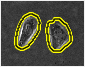}}}
         		\subfloat[]{
        \resizebox*{.45\linewidth}{!}{	\includegraphics[width=.5\textwidth]{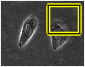}%    
         		\includegraphics[width=.5\textwidth]{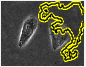}}}\\
         		\subfloat[]{
        \resizebox*{.45\linewidth}{!}{	\includegraphics[width=.5\textwidth]{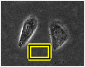}%    
         		\includegraphics[width=.5\textwidth]{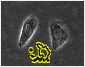}}}
         		\subfloat[]{\label{correct3}
        \resizebox*{.45\linewidth}{!}{	\includegraphics[width=.5\textwidth]{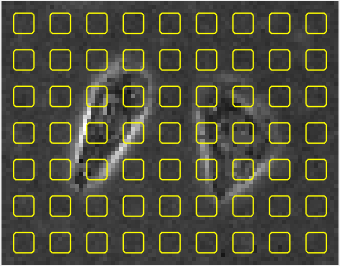}%    
         		\includegraphics[width=.5\textwidth]{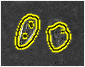}}}\\
        
        			 \caption[Sensitivity to initialization for the proposed variational model.]{Sensible initialization of the evolving front plays a crucial role in the segmentation results for the proposed model.  When initial contours encompass the target objects, as in \sref{correct1} and \sref{correct2}, meaningful results are obtained. In \sref{correct3} an almost correct segmentation is obtained if the detected region inside the cell on the left is neglected.}
         			\label{im:inicondexperiments}
        \end{figure}

  From a PDE-Variational point of view the problem of image segmentation  is ill-posed. Although objects in static images do not change (in some sense, they are unique) regardless of the location chosen for the initial evolving curves, two starting fronts  $ u_0, v_0,  u_0\neq v_0$, often render different solutions. Hence, there seems to be an unclear link between the steady state of evolution equations or (local) minima of energies and the particular image being segmented. Our proposed model does not escape this phenomenon. We tested our model against the `two cells' image using a number of different initializations keeping all other parameters constant.
% % % %\newpage
   \noindent \Cref{im:inicondexperiments} shows eight different initializations. It is found that the model tends to yield meaningful segmentations when single, regular initial curves surround the target objects. Also worth noting are initial curves shown in \Cref{correct3}. Contours of this form are desirable when detection of many objects is sought after. In this particular instance, not only the two cells were detected but also some structures within the left cell are suggested by the model as potential objects of interest.
   
  \subsubsection{\titlecap{Performance on multimodal images}}
\begin{figure}
\centering
\subfloat[]{
\resizebox*{.45\linewidth}{!}{
    \includegraphics[width=\textwidth]{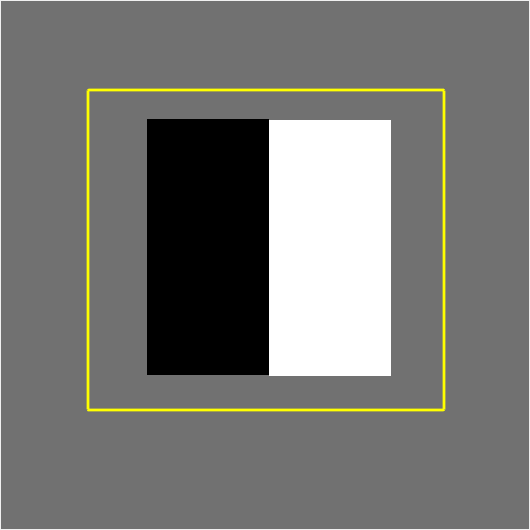}}}
\subfloat[]{\label{}
\resizebox*{.45\linewidth}{!}{
    \includegraphics[width=\textwidth]{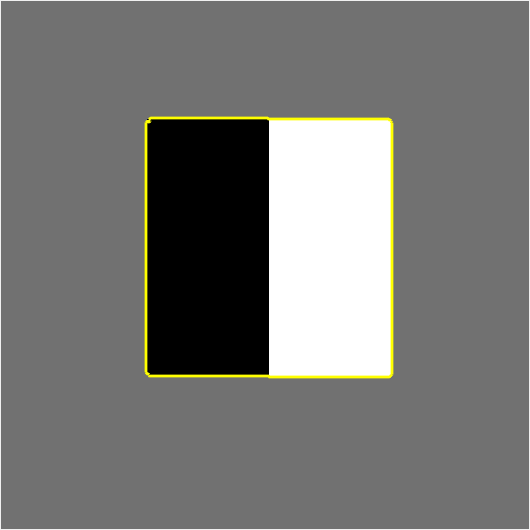}}} 
			\caption[Segmentation of trimodal image with triple junction.]{Segmentation of trimodal image with triple junction. (a) Original image and initial curve; (b) result by proposed model. Boundaries of the outer square are detected.}
	\label{im:trimodal}
\end{figure}

\begin{figure}
\centering
\subfloat[]{\label{fig:black2}
\resizebox*{.9\linewidth}{!}{
    \includegraphics[width=.45\textwidth]{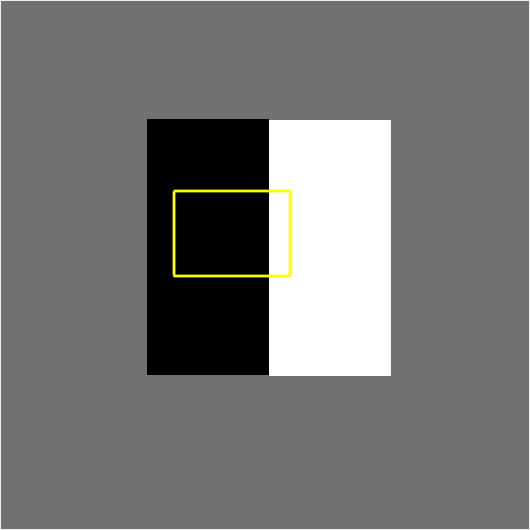}\includegraphics[width=.45\textwidth]{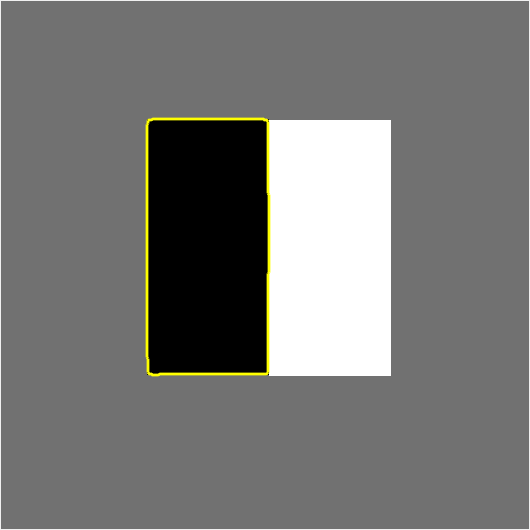}}}\\
    \subfloat[]{\label{fig:white2}
\resizebox*{.9\linewidth}{!}{
    \includegraphics[width=.45\textwidth]{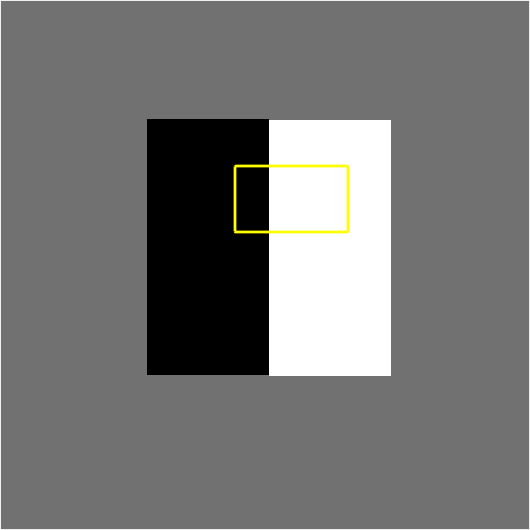}\includegraphics[width=.45\textwidth]{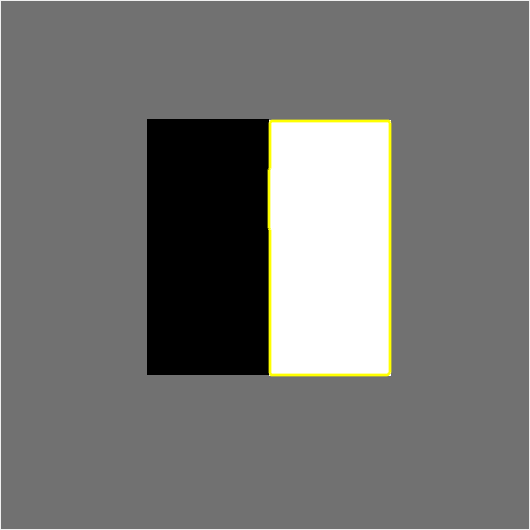}}} 
    \caption[Selective segmentation in trimodal image with triple junction.]{Selective segmentation in trimodal image with triple junction. Individual detection of \sref{fig:black2} black and \sref{fig:white2} white rectangle.}
 	\label{im:trimodal1}
\end{figure}

  Images containing more than one region of interest with different color/graylevel intensity are referred to as \textit{multimodal} images.  In this section we test our variational model on these images.
  \Cref{im:trimodal} shows our first experiment. The image consists of a square split into two black and white rectangles and the gray background. The three regions are adjacent to each other, a configuration commonly referred to as a \textit{triple junction}. One could be interested in separating the square from the background or either rectangle. Starting from a front surrounding the target objects, our model is able to detect the boundary of the square. This is possible due to incorporation of edge information. Use of region statistics alone fails as the gray level averages inside and outside the initial front are approximately equal.  
 The model is able to detect either rectangle as well, starting from an overlapping front mostly contained within the target rectangle (\Cref{im:trimodal1}). 

\begin{figure}[!htbp]
\centering
\subfloat[]{\label{fig:wh11f}
\resizebox*{.9\linewidth}{!}{
    	\includegraphics[width=.45\textwidth]{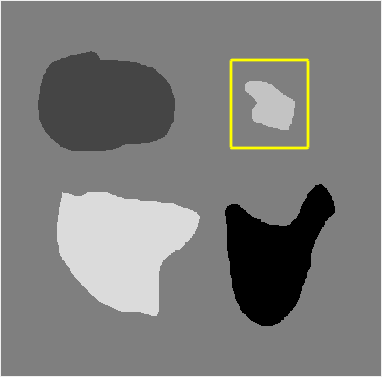}
    	\includegraphics[width=.45\textwidth]{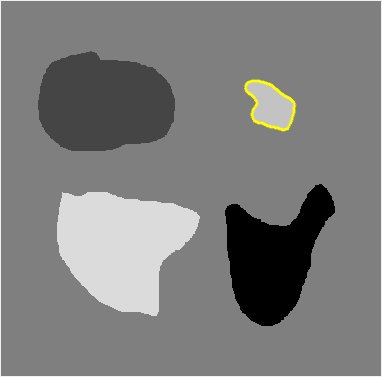}}}\\
    \subfloat[]{\label{fig:white32}
\resizebox*{.9\linewidth}{!}{
    	\includegraphics[width=.45\textwidth]{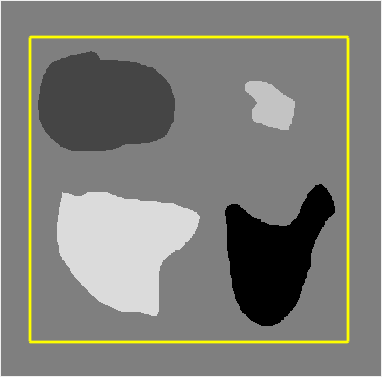}
    	\includegraphics[width=.45\textwidth]{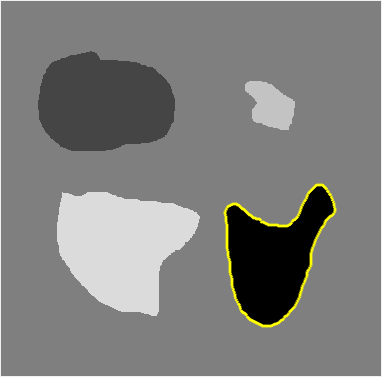}}} 
		\caption[Segmentation of image with 4 regions]{\sref{fig:wh11f} Selective detection is still possible, \sref{fig:white32} but detection of all objects in the scene is not.}
	\label{im:trimodal2}
\end{figure}

  Even if the objects of interest are not adjacent as in the example above, our model is still able to separate each object from the background, although selectively, as \Cref{im:trimodal2} illustrates. The model might fail if no prior information is available if a number objects are to be detected simultaneously with one single front. 
 \subsection{Sensitivity to Noise}
\begin{figure}[!htbp]
%file: noisy.tex
\centering
\resizebox*{.9\linewidth}{!}{
	\includegraphics[width=.45\textwidth]{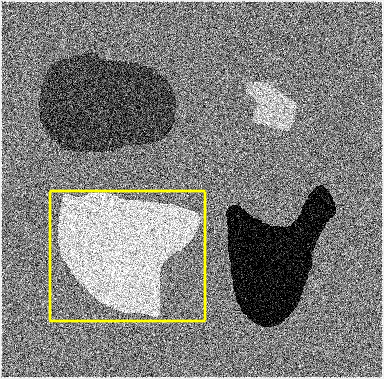} \includegraphics[width=.45\textwidth]{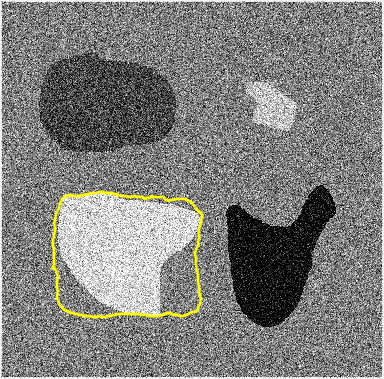}}
		\caption[Effect of noise.]{Effect of noise. Noise removal techniques ought to be implemented before segmentation.}
		\label{fig:noise}
\end{figure}

 We have shown that our segmentation model enjoys great versatility in terms of both the class of images it can process and its scope. Its main limitation is perhaps its sensitivity to noisy features in images. Being a hybrid model, spurious objects can significantly hamper its effectiveness as such anomalies are interpreted as possible boundaries of an object by the edge function $ g $, \Cref{fig:noise}. Therefore, preprocessing of noisy images using noise removal techniques \citep{Rudin92,Perona94} is recommended before the implementation of our variational model. 
 
\section*{Funding}
 This work was partially supported by FONDOCYT (Dominican Republic) under Grant 2016-2017-160017-2018.
\clearpage
\bibliographystyle{plainnat}
\bibliography{main}

\begin{thebibliography}{19}
\providecommand{\natexlab}[1]{#1}
\providecommand{\url}[1]{\texttt{#1}}
\expandafter\ifx\csname urlstyle\endcsname\relax
  \providecommand{\doi}[1]{doi: #1}\else
  \providecommand{\doi}{doi: \begingroup \urlstyle{rm}\Url}\fi

\bibitem[Adalsteinsson and Sethian(1995)]{Adalsteinsson95}
David Adalsteinsson and James~A Sethian.
\newblock A fast level set method for propagating interfaces.
\newblock \emph{Journal of computational physics}, 118\penalty0 (2):\penalty0
  269--277, 1995.

\bibitem[Aubert et~al.(2003)Aubert, Michel, Faugeras, and
  Jehan-Besson]{Aubert03}
Gilles Aubert, Barlaud Michel, Olivier Faugeras, and St\'ephanie Jehan-Besson.
\newblock Image segmentation using active contours: Calculus of variations or
  shape gradients?
\newblock \emph{SIAM Journal on Applied Mathematics}, 63\penalty0 (6):\penalty0
  2128--2154, 2003.
\newblock ISSN 00361399.
\newblock URL \url{http://www.jstor.org/stable/4096078}.

\bibitem[De~Vaucouleurs(1959)]{DeVaucouleurs59}
G\'erard De~Vaucouleurs.
\newblock Classification and morphology of external galaxies.
\newblock In \emph{Astrophysik IV: Sternsysteme/Astrophysics IV: Stellar
  Systems}, pages 275--310. Springer, 1959.

\bibitem[Demirkaya et~al.(2009)Demirkaya, Asyali, and Sahoo]{Demirkaya09}
Omer Demirkaya, Musa~Hakan Asyali, and Prasanna Sahoo.
\newblock \emph{Image processing with MATLAB : applications in medicine and
  biology}.
\newblock CRC Press, Boca Raton, 2009.
\newblock ISBN 9780849392467 0849392462.

\bibitem[Evans and Spruck(1991)]{Evans91}
Lawrence~C. Evans and Joel Spruck.
\newblock Motion of level sets by mean curvature {I}.
\newblock \emph{J. Diff. Geom}, 33\penalty0 (3):\penalty0 635--681, 1991.

\bibitem[Farag(2014)]{Farag14}
Aly~A. Farag.
\newblock \emph{BIOMEDICAL IMAGE ANALYSIS: STATISTICAL AND VARIATIONAL
  METHODS}.
\newblock Cambridge Univ Press, Cambridge, 2014.

\bibitem[Gonzalez et~al.(2004)Gonzalez, Woods, and Eddins]{Gonzalez04}
Rafael~C. Gonzalez, Richard~E. Woods, and Steven~L. Eddins.
\newblock \emph{Digital Image processing using MATLAB}.
\newblock Pearson Prentice Hall, Upper Saddle River, N. J., 2004.
\newblock ISBN 0130085197.

\bibitem[Kass et~al.(1988)Kass, Witkin, and Terzopoulos]{Kass88}
Michael Kass, Andrew Witkin, and Demetri Terzopoulos.
\newblock Snakes: Active contour models.
\newblock \emph{International journal of computer vision}, 1\penalty0
  (4):\penalty0 321--331, 1988.
\newblock ISSN 0920-5691.

\bibitem[Lankton and Tannenbaum(2008)]{Lankton08}
Shawn Lankton and Allen Tannenbaum.
\newblock Localizing region-based active contours.
\newblock \emph{Image Processing, IEEE Transactions on}, 17\penalty0
  (11):\penalty0 2029--2039, 2008.
\newblock ISSN 1057-7149.

\bibitem[Marques(2011)]{Marques11}
Oge Marques.
\newblock \emph{Practical image and video processing using MATLAB}.
\newblock Wiley-IEEE Press, Hoboken, NJ, 2011.
\newblock ISBN 9781118093467 1118093461 9781118093481 1118093488 9781118093474
  111809347X.
\newblock URL \url{http://site.ebrary.com/id/10494631}.

\bibitem[Mumford and Shah(1989)]{Shah89}
David Mumford and Jayant Shah.
\newblock Optimal approximations by piecewise smooth functions and associated
  variational problems.
\newblock \emph{Communications on pure and applied mathematics}, 42\penalty0
  (5):\penalty0 577--685, 1989.

\bibitem[Osher and Sethian(1988)]{Osher88}
Stanley Osher and James~A. Sethian.
\newblock Fronts propagating with curvature-dependent speed: algorithms based
  on hamilton-jacobi formulations.
\newblock \emph{Journal of computational physics}, 79\penalty0 (1):\penalty0
  12--49, 1988.
\newblock ISSN 0021-9991.

\bibitem[Perona et~al.(1994)Perona, Shiota, and Malik]{Perona94}
Pietro Perona, Takahiro Shiota, and Jitendra Malik.
\newblock Anisotropic diffusion.
\newblock In \emph{Geometry-driven diffusion in computer vision}, pages 73--92.
  Springer, 1994.

\bibitem[Rudin et~al.(1992)Rudin, Osher, and Fatemi]{Rudin92}
Leonid~I Rudin, Stanley Osher, and Emad Fatemi.
\newblock Nonlinear total variation based noise removal algorithms.
\newblock \emph{Physica D: Nonlinear Phenomena}, 60\penalty0 (1):\penalty0
  259--268, 1992.

\bibitem[Sethian(1999)]{Sethian99}
James~A. Sethian.
\newblock \emph{Level set methods and fast marching methods: evolving
  interfaces in computational geometry, fluid mechanics, computer vision, and
  materials science}, volume~3.
\newblock Cambridge university press, 1999.

\bibitem[Sethian(2001)]{Sethian01}
James~A Sethian.
\newblock Evolution, implementation, and application of level set and fast
  marching methods for advancing fronts.
\newblock \emph{Journal of Computational Physics}, 169\penalty0 (2):\penalty0
  503--555, 2001.

\bibitem[Tikhonov and Arsenin(1977)]{Tikhonov77}
Andrey~Nikolaevich Tikhonov and Vasiliy~Yakovlevich Arsenin.
\newblock \emph{Solutions of ill-posed problems}.
\newblock Vh Winston, 1977.

\bibitem[Yezzi et~al.(1999)Yezzi, Tsai, and Willsky]{Yezzi}
A.~Yezzi, A.~Tsai, and A.~Willsky.
\newblock A statistical approach to snakes for bimodal and trimodal imagery.
\newblock In \emph{Computer Vision, 1999. The Proceedings of the Seventh IEEE
  International Conference on}, volume~2, pages 898--903 vol.2, 1999.
\newblock \doi{10.1109/ICCV.1999.790317}.

\bibitem[Yezzi et~al.(2002)Yezzi, Tsai, and Willsky]{Yezzi02}
Anthony Yezzi, Andy Tsai, and Alan Willsky.
\newblock A fully global approach to image segmentation via coupled curve
  evolution equations.
\newblock \emph{Journal of Visual Communication and Image Representation},
  13\penalty0 (1):\penalty0 195--216, 2002.

\end{thebibliography}

\end{document}